# Study on S2 Flow Path Design and Three-dimensional Numerical Simulation Parameter Calibration in Axial Compressor

Yang Aobo [1*], Yan An [2**], Chen Jiang [3]

[1] School of Energy and Power Engineering, Beihang University, Beijing, China
[2] School of Energy and Power Engineering, Beihang University, Beijing, China
[3] School of Energy and Power Engineering, Beihang University, Beijing, China

**Abstract**

Aerodynamic design process of multi - stage axial flow compressor usually uses the way that combines the S2 flow design and three-dimensional CFD numerical simulation analysis. Based on Mr. Wu Zhonghua's " Three-dimensional flow theory ", aiming at the S2 flow design matching parameters and the three-dimensional CFD numerical simulation data, through autonomous programming, the S2 design parameter distribution and the corresponding CGNS format CFD calculation results are extracted. Then make the comparative analysis of the two and provide the modification suggestion of the design. The examples have been tested by the comparison and correction of the eight-stage axial flow compressor calculation and finally improve the design performance of the compressor design. The design adiabatic efficiency increases by 0.5%. The surge margin increases by more than 5%. The validity and feasibility of the method are verified.

## 1. Introduction

Currently multistage axial-flow compressors usually use general aerodynamic design methods. According to "Three-dimensional flow theory" proposed by Wu Zhonghua[1,2], the process is as follows: First, carry out S2 flow path design of compressor, conduct blade modeling according to the result of S2 flow path design. Then, check and analyze the three-dimensional numerical simulations. If the result doesn't meet the expected requirements, modify the blade design and carry out the verification of full three-dimensional numerical simulations again. Therefore 3D numerical simulation result and S2 design parameters iterative correction has become an important procedure of the design for compressors. By comparing the differences between 3D numerical simulating result and S2 design, we can modify the blade modeling and achieve the expected working performance of compressors in order to improve the design accuracy.

In recent years, owing to the rapid development of computer technology and significant progress of computational fluid dynamics(CFD) theories, CFD applied in turbomachinery is also developing rapidly and has been fully involved in the design and developing process of compressors, reducing developing costs and cycle. The main CFD programs that are active in domestic and foreign academic and engineering field include: Denton[3] program, ADPAC program, APNASA program from NASA[4], programs developed by Tsinghua University, Chinese Academy of Sciences, Beihang University and other research institutes[5], and some mature business software, such as CFX, NUMECA, Fluent, etc. All of them have been applied in some practical research and projects and have got lots of achievements.

In this paper, the research object is an eight-stage axial flow compressor. By analyzing the 3D numerical simulating CGNS format file of eight-stage axial compressor, the extracting method of multistage axial flow compressor full 3D numerical simulating data and data processing method are proposed, the comparing and analyzing result of S2 flow design and 3D numerical simulating result are realized, facilitating the designers to find differences, modify the blade modeling design, iteratively check the 3D calculation and S2 design, and finally improve the performance of axial flow compressors.

## 2. Extraction Method of CFD calculating data and S2 Flow Path Design Parameters

Since there are various data types and large amount of the eight-stage axial compressor 3D simulation results, we choose to extract the three-dimensional simulation data by analyzing the three-dimensional numerical simulation of the compressor CGNS standard format. And the simulation data are compared with the extracted S2 flow path design parameters. Hao Jia[6], from Beijing Jiaotong University, has researched aviation CFD calculating data automatic processing technology based on CGNS. It is used to compare the data post-processing calculation in the design of aircrafts and the data in wind tunnel experiment. After multiple tests and practices, the effectiveness and efficiency of data extraction processing are verified.

* Presenting Author: Yang Aobo, yangaobo1995@163.com, Student, BUAA.
** Corresponding Author: Chen Jiang, chenjiang27@buaa.edu.cn, Professor, BUAA.

## 2.1. CGNS Standard Storage Format

CGNS standard storage format[7] was created by the Boeing Company and NASA in 1994, trying to standardize the input and output of CFD, including grid (structured and unstructured grid), flow solutions, connectivity, BC and auxiliary information. The purpose is to eliminate the translating work between the machines and CFD programs. CGNS uses the ADF (Advanced Data Format) system, which creates a binary file, making it more convenient to cross the computer platform, and makes CGNS become the standardized storage format file of impeller mechanical full 3D CFD numerical simulation result.

## 2.2. Extraction Method of CGNS Format Data Based on CFD Numerical Simulation

Establish the required function library and extraction environment for CGNS extraction through Fortran language based on the Linux system[8,9]. The extraction process is divided into four steps: First, realize the compilation of HDF5 library and CGNS function library. Second, create a 3D grid file, and achieve the correspondence of CGNS file data and grid. Third, read the grid file and generate the required mod file. Fourth, realize the extraction of CGNS data. Fig.1 is the flow chart of data extraction of S2 flow path design and 3D numerial simulation.

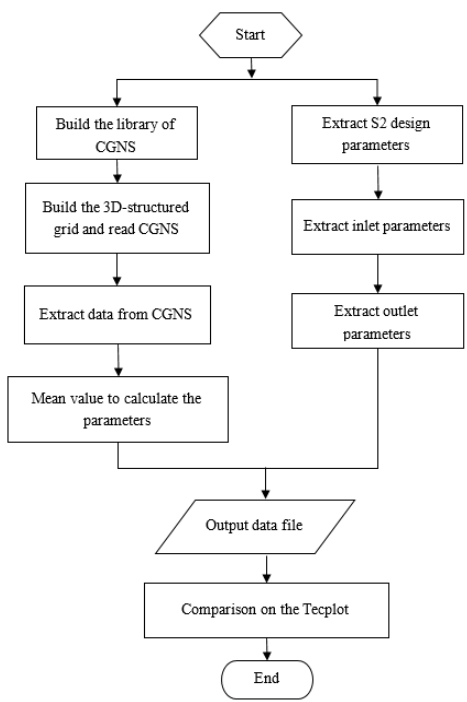

Fig.1. Data extraction process of S2 flow path design and 3D numerical simulation

After all the three-dimensional numerical simulation data of CGNS are extracted, it is necessary to carry out the axial average of the blades at the position of inlet and outlet for a specific physical quantity, in order to ensure that the three-dimensional numerical simulation results after the averaging algorithm can be compared with the S2 flow path design parameters in the same direction of blade height. The data extraction and average algorithm are as follows:

First, read all the variable names in CGNS file and output the required statement meeting the Tecplot reading file. Second, according to the CGNS storage standard, read each domain and grid points of each base in turn, calculate the number of grid points of every domain, filter required data of meridional flow surface and deposit all the data of required physical quantity into two-dimensional arrays. Finally, find the grid point coordinates meeting S2 reference condition. The condition lies in the height of required leaf which meets the error range (R-ε<X<R+ε). Calculate the average value of each domain and output the data.

## 2.3. Extraction of S2 Flow Path Design Parameters

First, we need to get the title of axial compressor stages, calculate the leaf sorting number of corresponding CGNS format according to the quantitative relationship between every stator and rotor blade of compressors. Since the import and export data value of rotor and stator are different, they should be separated when being extracted. Then read required variables.

By completing the process of S2 flow path design parameters and three-dimensional numerical simulation parameters, we can find that CGNS data extraction is more accurate, detailed and could extract more data for the cross-section selection and grid reading than other three-dimensional numerical simulation data extraction methods. Also, the expression of parameter relationship is more complete and the CGNS data extraction has a wider range of application. Therefore, in terms of the comparison of S2 design and the three-dimensional numerical simulation result, especially the design and verification work of multi-stage axial flow compressors, which has a large amount of data and data types, the extraction of CGNS file data could provide a better result: the amount is more complete and the image is more detailed and accurate, making it better to give a reference for modifying the blade design of compressors.

## 2.4. Comparing Method of CFD Calculation and S2 Flow Path Design Parameters

In this paper, the researching object is an eight-stage axial compressor. After the extraction of S2 flow path design calculating data and the three-dimensional numerical simulation data, we need to compare the two types of data in the inlet and outlet of each stage in axial compressor on the same blade height. Through the illustration of Tecplot charts, the differences of the two could be visually analyzed and modifying suggestions are proposed. According to the suggestions, S2 flow design and the blade design are improved, achieving the anticipative designing requirements of the compressor and improving the performance of axial compressor design.

The comparing result of S2 flow path design and the three-dimensional numerical calculation is shown through Tecplot. The results are as follows: the designers can choose any axial compressor and any inlet or outlet domain of the stage. They can choose 3D calculation and S2 design parameters, of which the same physical quantity is required to be compared. The contrast images of the three-dimensional numerical simulation and S2 design result in the same blade length direction are shown. Users can observe the image and compare on the whole, and can also choose the data of specific points. In order to get a better view of the images, the blade length is nondimensionalized and is presented in the form of the percentage of the blade length.

## 3. S2 Flow Path Design and 3D Numerical Calculating Parameters Calibration of Eight-stage Axial Compressor

As CFD numerical simulation analysis says: S2 flow path design has basically achieved the designing targets. The next step lies in the modification and improvement of the leaf model.

Due to the excessive number of the blades in the eight-stage axial compressors, and there are too many reference physical quantities, only researching methods for parameters comparison and revision of S2 path flow design and 3D numerical analysis are illustrated in this paper. Some modifying suggestions are also proposed. Since axial compressors design calculating revision is the process of iterative comparison, the compressor shape is supposed to be modified repeatedly according to experience, therefore the specific process isn't described in detail in this paper.

### 3.1. Introduction of research objects

The entrance of the eight-stage axial compressor researched in this paper is under the standard atmospheric conditions, of which the total pressure is 101325Pa, the total temperature is 288℃, with intake in the axial direction. The original backpressure is 600000Pa, in the following step the backpressure will be enhanced to calculate. Fig.2 shows the eight-stage axial flow compressor blade design. Fig.3 shows the whole eight-stage flow compressor.

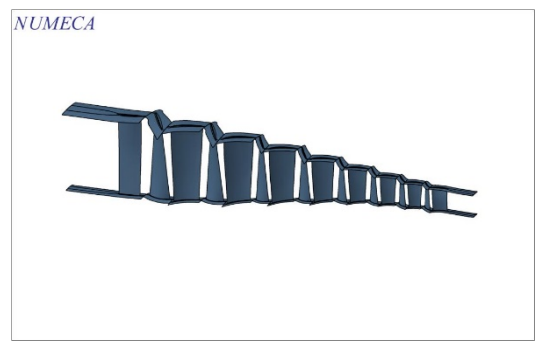

Fig.2.Design of the eight-stage axial-compressor blade row

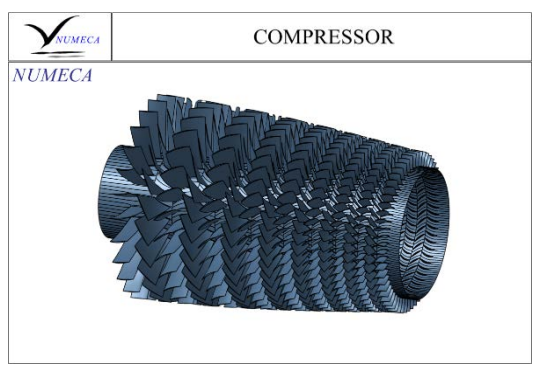

Fig.3.Eight-stage axial-compressor machine display

3D numerical simulation uses the Fine Turbo module of the software NUMECA. The fluid filling governing equations are 3D steady Reynolds average N-S equations, and S-A equation turbulence model is adopted. Central difference format is adopted for spatial spreading. The time derivative term is iteratively solved by the fourth-order Runge-Kutta method. There are 61 stations for dynamic blade of S2 path flow design, 17 for radial clearance and 49 for stationary blades.Fig.4.shows the NUMECA 3D calculating grids.

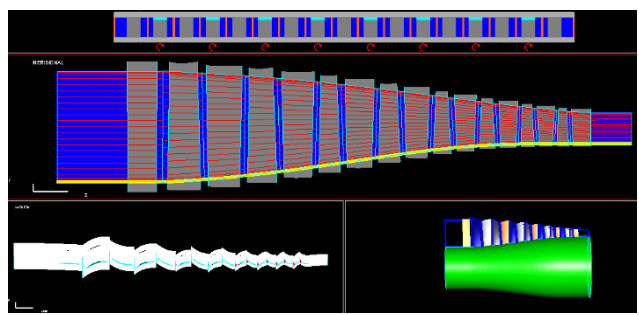

Fig.4.NUMECA 3D calculation

### 3.2. S2 Path Flow Design and 3D Numerical Calculating Parameters Calibration

By observing the comparing results of the eight-stage axial-flow compressor S2 flow design and the 3D numerical simulation data presented by Tecplot and analyzing the reason for difference as well as the lack of design, adjust the geometrical parameters such as compressor leaf type according to experience and redesign S2 path flow. Compare and analyze the new S2 path flow design and its 3D numerical simulating result. Repeat the above process iteratively until the S2 path flow design and 3D numerical calculation could achieve the maximum fitting, therefore the compressor can achieve the expected performance. After the repeatedly and iteratively modification for S2 design, the designing performance of compressors is finally improved, and the effectiveness and feasibility of the researching method S2 path flow design and 3D numerical analyzing parameter modification proposed in this paper.[10] Fig.4 illustrates the axial compressor aerodynamic design and the procedure of iterative modification.

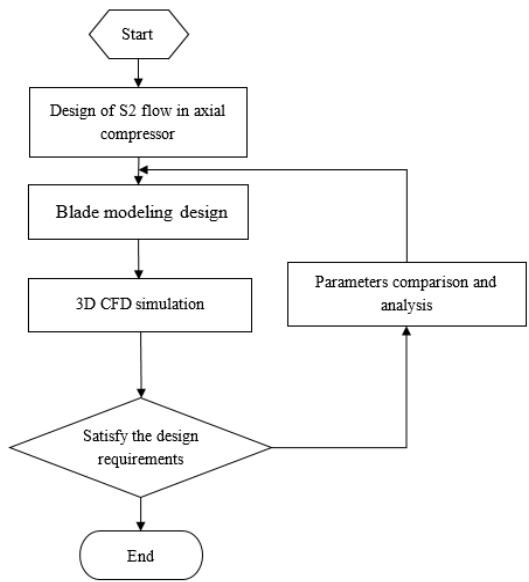

Fig.5.Aerodynamic design flow of axial compressor

By comparing each stage design pressure ratio and efficiency between 8-stage axial compressor and their 3D numerical simulation results, according to the Table 1(show parts because of space). We can draw a conclusion that the design numerical simulation result basically achieve the design goal of each stage. If we want to continue to elevate the performance of compressor, we should concentrate more on the NO.1 and NO.8 of the compressor, especially the No.8 stage.

Table 1.S2 and 3D stage performance comparasion

| Stage | 1 | 4 | 8 |
|---|---|---|---|
| S2 pressure ratio | 1.37 | 1.30 | 1.16 |
| 3D pressure ratio | 1.349 | 1.30 | 1.1544 |
| S2 efficiency | 0.9424 | 0.9388 | 0.9140 |
| 3D efficiency | 0.932 | 0.9451 | 0.9022 |

### 3.3. Suggestions of modification

(1) The relative mach distributed in the blade tip

The follows are the relative mach numbers distributed in the blade root(5% blade height), medium(50% blade height), tip(95% blade height), shown in fig.5,6, we can see that blade back of rotor in stage 1 exits supersonic shock(Mach is approximately larger than 1.3). We need to design in-depth and accurately in order to decrease the relative mach number. We may consider the subsonic aerodynamic design to improve the efficiency in order to decrease the loss of shock wave.

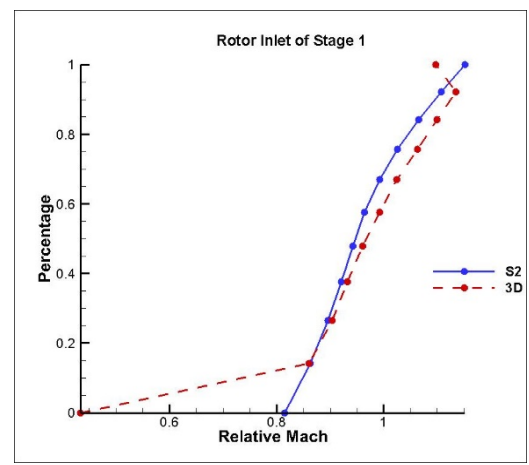

Fig.6.Comparasion between S2 and 3D of rotor inlet of stage 1

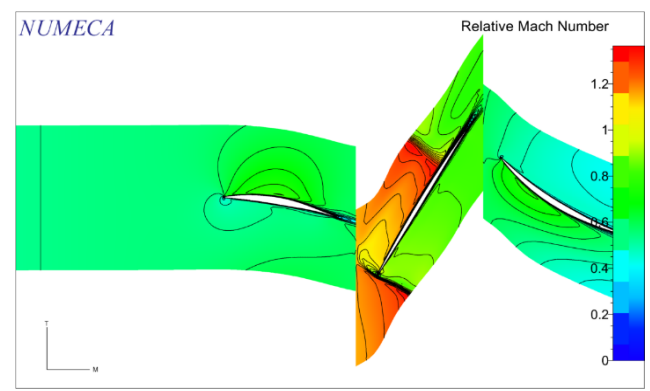

Fig.7.Relative mach number distribution in blade tip of stage 1 rotor

(2) Total pressure of comparasion

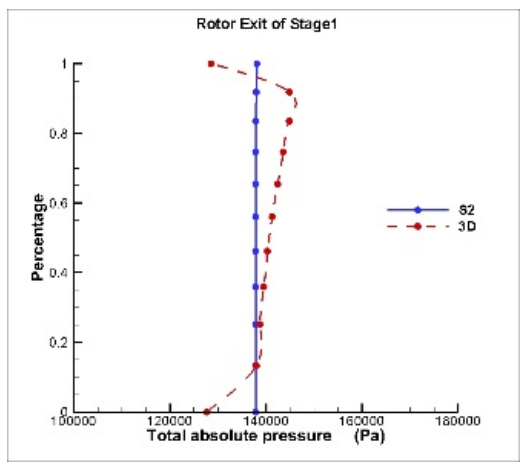

Fig.8.Total Pressure comparation of Rotor exit of Stage 1

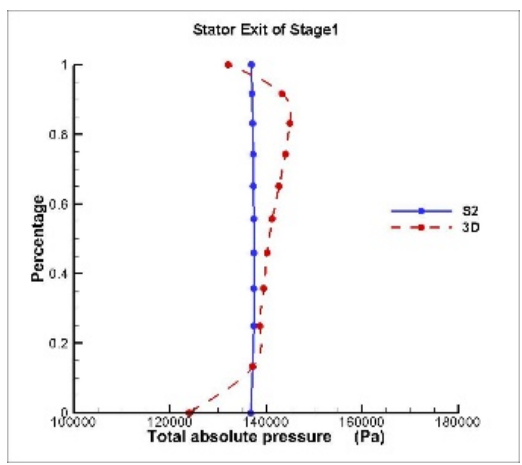

Fig.9.Total Pressure comparation of Stator exit of Stage 1

In the comparasion of inlet and outlet total pressure of stage1, the calculating result of S2 flow path design along the height of the blade is approxiamately 1389000Pa, lower than 3D result at the tip and medium of the blade, we need to adjust the blade attack angle and bowed angle to change the ratio of pressure increase.

The comparasion between S2 path flow paeameters and three-dimensional anaylsis data can be a useful auxiliary method that help designers to analyze the data and images together with three-dimensional calculation images. The above shows the way it works. And the blade modeling need iteractive modification according to the comparasion results of three-dimensional numerical analysis and the S2 flow path design.

## 4. Conclusion

The paper takes an eight-stage axial compressor as an example, achieve the objective of extraction and comparasion of S2 flow path design parameters and 3D numerical simulation data results and finish the study of parameters calibration of the two above. Finally, the axial compressor performance has been elevated, the design conversion efficiency increase 5%, the surge margin increased by more than 5%. The conclusions are as follows:

(1): The method of extracting data from CGNS standard format in the process of multi-stage axial compressor aerodynamic design with the advantages of large data amount and high calculating accuracy. The program developed applies to the extraction of all multi-stage axial compressors 3D numerical simulation CGNS format;

(2): Overall analysis conclude that the efficiency of 8-stage axial compressor can be improved by refining the design,we need to concentrate on the improvement of the first and the last stage of compressor and improve the margin in low speed undesigned conditions.